\documentclass[11pt]{amsart}

\usepackage{amsmath,amssymb,amsthm}
\usepackage[hidelinks]{hyperref}
\usepackage{booktabs}
\usepackage{url}

\newtheorem{theorem}{Theorem}[section]
\newtheorem{proposition}[theorem]{Proposition}
\newtheorem{lemma}[theorem]{Lemma}
\newtheorem{corollary}[theorem]{Corollary}
\theoremstyle{definition}
\newtheorem{definition}[theorem]{Definition}
\theoremstyle{remark}
\newtheorem{remark}[theorem]{Remark}

\newcommand{\Q}{\mathbb{Q}}
\newcommand{\Z}{\mathbb{Z}}
\newcommand{\C}{\mathbb{C}}
\newcommand{\R}{\mathbb{R}}
\DeclareMathOperator{\cl}{cl}
\DeclareMathOperator{\Imz}{Im}
\newcommand{\diam}{\diamond}
\title[A transmitted form of Bloch's diamond formula]{A normalization
issue in a transmitted form of Bloch's diamond formula:\\
anti-invariant subgroup versus regulator quotient}

\author{Huimin Zheng}
\address{College of Information and Network Engineering,
Anhui Science and Technology University,
Fengyang, Anhui 233100, P.~R.~China}
\email{zhenghm@ahstu.edu.cn}

\subjclass[2020]{Primary 11G55, 11G40; Secondary 11G05, 11R06}
\keywords{Bloch's diamond formula, elliptic dilogarithm, Beilinson
regulator, Mahler measure, anti-invariant homology, high-precision
numerical verification}

\date{November 2, 2026}

\begin{document}

\begin{abstract}
Bloch's regulator formula for an elliptic curve is often quoted with a generator of the anti-invariant subgroup $H_1(E,\mathbb{Z})^-$ and coefficient $1$ in front of the elliptic dilogarithm of the diamond product. We compare this wording with Brunault's published intersection-product formula. The natural map from the anti-invariant subgroup to the regulator quotient $H_1(E,\mathbb{Z})/H_1(E,\mathbb{Z})^+$ has index $1$ when $E(\mathbb{R})$ has two components and index $2$ when it has one. Thus, up to orientation, the coefficient is $1$ for a quotient generator and $\iota(E)\in\{1,2\}$ for a subgroup generator. This is an elementary consequence of the published formula, not a new regulator theorem. We illustrate the distinction on 11.a3 using an externally proved regulator anchor and the Brunault--Bertin special values, and on 20.a4 by high-precision numerical integration. We also correct our earlier quotation of Touafek: his printed identity is $\pi r=D_E((f)\diamond(g))$, with $r=(2\pi)^{-1}\int_{\gamma^-}\eta$. It is compatible with the subgroup normalization when $\Delta<0$; when $\Delta>0$ the corresponding coefficient is $1/2$, up to orientation. We do not reassess his isogeny theorem or his final exotic relation here. No new proof of Boyd's conductor-11 split-integral identity is included.
\end{abstract}

\maketitle

%=====================================================================
\section*{Declaration on the use of AI tools}
The author used \emph{Kimi} (Moonshot AI) for computational exploration
and preparation of the earlier manuscript, and \emph{Codex} (OpenAI)
to assist with checking quotations, revising the text, and preparing
this corrected version. The author is responsible for the mathematical
content, the references, and the interpretation of all computations.
The numerical script and its run log accompany the source distribution
(see Section~\ref{sec:numerics}).

%=====================================================================
\section{Introduction}

Let $E/\Q$ be an elliptic curve and let $f,g\in\Q(E)^\times$ be such
that the symbol $\{f,g\}$ lies in $K_2(E)\otimes\Q$. Two real numbers
are attached to these data. The first is the \emph{regulator integral}
\begin{equation}\label{eq:regint}
  r(\{f,g\})[\gamma]\;=\;\int_\gamma\eta(f,g),\qquad
  \eta(f,g)=\log|f|\,d\arg g-\log|g|\,d\arg f,
\end{equation}
over a loop $\gamma$ on $E(\C)$. The second is the \emph{diamond value}
\begin{equation}\label{eq:diamondvalue}
  D_E\bigl((f)\diam(g)\bigr),\qquad
  (f)\diam(g)=\textstyle\sum_{i,j}m_i n_j\,\cl(S_i-T_j),
\end{equation}
where $(f)=\sum_i m_i(S_i)$, $(g)=\sum_j n_j(T_j)$, and $D_E$ is
Bloch's elliptic dilogarithm (Definition~\ref{def:DE} below). That the
two agree, for a suitable loop $\gamma$, is a theorem attributed to
Bloch \cite{Bloch}; in the form in which it is
routinely cited in the Mahler-measure literature initiated by Boyd's
tables \cite{Boyd98}, for instance in
\cite[Thm.~6]{LalinRam}, the identity reads
\begin{equation}\label{eq:transmitted}
  r(\{f,g\})[\gamma]\;=\;D_E\bigl((f)\diam(g)\bigr),
  \qquad\text{$\gamma$ a generator of $H_1(E,\Z)^-$,}
\end{equation}
where $H_1(E,\Z)^-$ is the anti-invariant subgroup, the $(-1)$-eigenspace
of complex conjugation on $H_1(E(\C),\Z)$.

In the course of the certified computations accompanying \cite{Paper1}
we observed that \eqref{eq:transmitted} fails on the curve 11.a3 by a
factor of exactly $2$: both sides of \eqref{eq:transmitted} are
\emph{proved} quantities there (the left side modulo the anchor theorem
of \cite{Paper1}; the right side by Brunault's and Bertin's
evaluations), and they differ by $2$. The purpose of this note is to
pin down the source of the discrepancy, to state and prove the correct
general form of the identity, and to document the phenomenon.

The source is a lattice, not a sign or a convention for $D_E$. The
pairing \eqref{eq:regint} vanishes on the invariant subgroup
$H_1(E,\Z)^+$ (because $c^*\eta=-\eta$), so the regulator is naturally
a function on the \emph{regulator quotient}
$H_1(E,\Z)/H_1(E,\Z)^+$. The anti-invariant subgroup maps into this
quotient with index
\begin{equation}\label{eq:iota}
  \iota(E)\;=\;
  \begin{cases}
    1, & \text{$E(\R)$ has two components ($\Delta>0$),}\\[2pt]
    2, & \text{$E(\R)$ has one component ($\Delta<0$).}
  \end{cases}
\end{equation}
Our main result is that \eqref{eq:transmitted} holds with factor $1$
for a generator of the quotient (Theorem~\ref{thm:main}), and
consequently holds for the subgroup generator exactly as transmitted
when $\Delta>0$, but with the right-hand side multiplied by $2$ when
$\Delta<0$ (Corollary~\ref{cor:subgroup}). The proof is short:
Brunault's intersection-product formula for the regulator integral
\cite[Prop.~19]{BrunaultSiegel}\footnote{Numbered Proposition~18 in
the arXiv version (arXiv:1504.08127v1); the published version and the
author's revised preprint of February~2016 number it Proposition~19,
and \S5.1 of the preprint became \S6.1.}
is a published, proved statement whose
normalization is intrinsic, and the index computation of
Section~\ref{sec:index} translates it into the subgroup/quotient
language.

\paragraph{Positioning.}
The intersection-product normalization, including the factor $2$ in
applications with $\Delta<0$, is already published
\cite[Prop.~19 and \S6.1]{BrunaultSiegel}. We do not claim a new
regulator theorem. Our purpose is to make explicit the lattice
distinction needed when interpreting the subgroup wording of
\eqref{eq:transmitted}. With compatible orientations, that wording
has the wrong magnitude when $\Delta<0$ and the diamond value is
nonzero. The proof below is an elementary consequence of the published
formula.

The examples discussed here distinguish three situations.
First, the conductor-$17$ calculation of \cite{LalinRam} takes place
on a curve with $\Delta>0$, where the two lattices coincide.
Second, identities stated purely in terms of $D_E$-values, such as
Bertin's proved exotic relations \cite{Bertin,BertinCrelle}, do not
themselves specify an integral cycle; the lattice correction alone
does not contradict them.
Third, Touafek uses a different analytic normalization: his
Definition~2 includes $(2\pi)^{-1}$, and his Theorem~1 and
Proposition~1 print $\pi r$, not $r$, on the left.
On his curve 20.a4, which has $\Delta<0$, this is compatible with
Corollary~\ref{cor:subgroup}, up to orientation
(Remark~\ref{rem:touafek}).
The earlier draft of this note omitted that printed $\pi$ and drew
an incorrect conclusion about two independent inconsistencies.
That conclusion is withdrawn. The present note does not assess
Touafek's Theorem~2 or prove or disprove his final exotic relation.

\paragraph{Related work.}
The normalization theorem is independent of \cite{Paper1}: its proof
uses only Brunault's published formula and the elementary lattice
computation. The conductor-$11$ example cites the regulator-anchor
theorem of \cite{Paper1} as an external input. We do not give a proof
of Boyd's split-integral identity here, nor include the chain,
branch, intersection, or positivity arguments used for that identity.
A systematic audit of all related literature is outside the scope
of this note.

%=====================================================================
\section{Setup}\label{sec:setup}

Throughout, $E/\Q$ is an elliptic curve, $c$ denotes complex
conjugation acting on $E(\C)$, and
\[
  H_1(E,\Z)^\pm\;=\;\ker(c_*\mp1)\;\subset\;H_1(E(\C),\Z).
\]
We write $\eta(f,g)=\log|f|\,d\arg g-\log|g|\,d\arg f$, a real
$1$-form on $E(\C)$ minus the support of $(f),(g)$, satisfying
$c^*\eta(f,g)=-\eta(f,g)$. Consequently the regulator integral
\eqref{eq:regint} vanishes on $H_1(E,\Z)^+$: if $c(\gamma)=\gamma$
then $\int_\gamma\eta=\int_\gamma c^*\eta=-\int_\gamma\eta$
(cf.\ \cite[Rem.~4]{LalinRam}). It therefore descends to the
\emph{regulator quotient}
\[
  Q(E)\;=\;H_1(E,\Z)\big/\,H_1(E,\Z)^+,
\]
a free abelian group of rank $1$. (This is \emph{not} the coinvariant
group $H_G=H_1/(c_*-1)H_1$ of group homology: for $\Delta>0$ the
latter is $\Z a\oplus(\Z/2)[b]$ in the basis of
Lemma~\ref{lem:basis}. Precisely, writing $G=\langle c\rangle$ and
$\Z(1)$ for the sign representation, $Q(E)\cong
\bigl(H_1(E,\Z)\otimes\Z(1)\bigr)_G/\mathrm{torsion}$, the torsion-free
quotient of the sign-twisted coinvariants.) We emphasize the
distinction drawn in
\cite{Paper1}: $Q(E)$ and $H_1(E,\Z)^-$ are two lattices inside the
same rational line, and they differ by the index that is the subject of
this note.\footnote{Strictly speaking, the arithmetic regulator is
defined on $K_2(\mathcal{E})\otimes\Q$ for the N\'eron model
$\mathcal{E}$, a subgroup of $K_2(E)\otimes\Q$ cut out by finitely many
extra conditions \cite[Rem.~4]{LalinRam}. Our main theorem is proved
via \cite[Prop.~19]{BrunaultSiegel}, which applies to symbols in
$K_2(E)\otimes\Q$ itself, so nothing in this note depends on the
distinction.}

\begin{definition}[The elliptic dilogarithm]\label{def:DE}
Fix an orientation of $E(\R)$, let $\omega_1$ be the primitive real
period of the identity component in that direction, and choose
$\omega_2$ such that $\tau=\omega_2/\omega_1$ lies in the upper
half-plane with real part $0$ or $\tfrac12$ \textup{(}a period basis
compatible with the real structure\textup{)}; put
$q=\exp(2\pi i\tau)$ and uniformize $E(\C)\cong\C^\times/q^\Z$ by
$t\mapsto z=\exp(2\pi i t/\omega_1)$. Let
$D(z)=\Imz\bigl(\operatorname{Li}_2(z)\bigr)+\arg(1-z)\log|z|$ be the
Bloch--Wigner dilogarithm. For $P\in E(\C)$ corresponding to
$z\in\C^\times/q^\Z$,
\begin{equation}\label{eq:DE}
  D_E(P)\;=\;\sum_{n\in\Z}D(q^n z),
\end{equation}
extended by linearity to divisors; it factors through
$\Z[E(\C)]^-=\Z[E(\C)]/\langle(S)+(-S)\rangle$. With the uniformization
normalized as above --- as in \cite[Def.~5]{LalinRam} and
\cite[\S6]{BrunaultSiegel} --- the series depends only on the chosen
orientation of $E(\R)$, so $D_E$ is a well-defined real-valued
function (in particular, \emph{no factor of $2$ is hidden in its
definition}; the series \eqref{eq:DE} is the one used in all of our
computations).
\end{definition}

%=====================================================================
\section{The index computation}\label{sec:index}

\begin{lemma}[The real homology basis]\label{lem:basis}
Let $E/\R$ be an elliptic curve. There is a $\Z$-basis $(a,b)$ of
$H_1(E(\C),\Z)$ with intersection number $a\cdot b=1$ such that
\begin{enumerate}
  \item if $\Delta>0$ \textup{(}$E(\R)$ has two components\textup{)},
        $c_*a=a$ and $c_*b=-b$; hence $H_1(E,\Z)^+=\Z a$,
        $H_1(E,\Z)^-=\Z b$;
  \item if $\Delta<0$ \textup{(}$E(\R)$ has one component\textup{)},
        $c_*a=a$ and $c_*b=a-b$; hence $H_1(E,\Z)^+=\Z a$ and
        $H_1(E,\Z)^-=\Z\,(a-2b)$.
\end{enumerate}
\end{lemma}

\begin{proof}
Uniformize $E(\C)\cong\C/\Lambda$ with $c$ induced by $z\mapsto\bar z$
(the shape of the period lattice of a real elliptic curve is standard;
see e.g.\ \cite[Ch.~3]{Cremona}).
For $\Delta>0$ the lattice is rectangular, $\Lambda=\Z\omega_1+\Z\omega_2$
with $\omega_1>0$, $\omega_2\in i\R_{>0}$; take $a,b$ to be the classes of
$s\mapsto s\omega_1$, $s\mapsto s\omega_2$ ($0\le s\le1$). For $\Delta<0$
one may take $\Lambda=\Z\omega_1+\Z\omega_2$ with $\omega_1>0$ and
$\omega_2=\omega_1/2+it$ ($t>0$); then
$\overline{\omega_2}=\omega_1-\omega_2$, so $c_*b=a-b$. In case (2),
$\alpha a+\beta b$ is anti-invariant iff
$(\alpha+\beta)a-\beta b=-\alpha a-\beta b$, i.e.\ $2\alpha+\beta=0$,
giving $H_1(E,\Z)^-=\Z(a-2b)$. Case (1) is immediate.
\end{proof}

\begin{lemma}[The index]\label{lem:index}
The natural map $H_1(E,\Z)^-\to Q(E)$ \textup{(}injective, since
$H_1(E,\Z)^-\cap H_1(E,\Z)^+=0$\textup{)} has image of index
$\iota(E)$, with $\iota(E)$ as in \eqref{eq:iota}. More precisely, with
the basis of Lemma~\ref{lem:basis}, $Q(E)=\Z\,[b]$, and the class of
the anti-invariant generator satisfies
\[
  [b]=[b]\quad(\Delta>0),\qquad [a-2b]=-2\,[b]\quad(\Delta<0).
\]
\end{lemma}

\begin{proof}
$H_1(E,\Z)^+=\Z a$ in both cases, so $Q(E)=\Z[b]$. In case (1) the
anti-invariant generator $b$ projects to $[b]$; in case (2) the
anti-invariant generator $a-2b$ projects to $-2[b]$.
\end{proof}

\begin{lemma}[The pairing with the invariant class]\label{lem:pairing}
Let $E^+\in H_1(E,\Z)^+$ be a generator. The intersection pairing
$(\gamma,\bar\delta)\mapsto\gamma\cdot\tilde\delta$ \textup{(}any lift
$\tilde\delta$ of $\bar\delta\in Q(E)$\textup{)} is a well-defined
perfect pairing $H_1(E,\Z)^+\times Q(E)\to\Z$. Hence a generator
$\bar\gamma$ of $Q(E)$ has a lift $\gamma$ with
$E^+\cdot\gamma=\pm1$.
\end{lemma}

\begin{proof}
Well-definedness: if $\tilde\delta$ is replaced by
$\tilde\delta+\delta^+$ with $\delta^+\in H_1(E,\Z)^+$, the pairing
changes by $\gamma\cdot\delta^+$, which vanishes: complex conjugation
is anti-holomorphic, hence reverses the orientation of $E(\C)$ and
sends the intersection number to its negative, while $c_*\gamma=\gamma$
and $c_*\delta^+=\delta^+$, so
$\gamma\cdot\delta^+=(c_*\gamma)\cdot(c_*\delta^+)=-(\gamma\cdot\delta^+)$.
Perfectness: in the basis of Lemma~\ref{lem:basis}, $E^+=\pm a$ and
$a\cdot b=1$, so $a\cdot[b]=1$.
\end{proof}

%=====================================================================
\section{The corrected statement}\label{sec:corrected}

\begin{theorem}[The diamond formula on the regulator
quotient]\label{thm:main}
Let $E/\Q$ be an elliptic curve, $\{f,g\}\in K_2(E)\otimes\Q$, and
$\bar\gamma$ a generator of $Q(E)=H_1(E,\Z)/H_1(E,\Z)^+$. Then
\begin{equation}\label{eq:main}
  r(\{f,g\})[\bar\gamma]=\int_{\bar\gamma}\eta(f,g)
  =\pm\,D_E\bigl((f)\diam(g)\bigr),
\end{equation}
where the sign corresponds to the choice of generator of $Q(E)$: once
the orientation of $E(\R)$ defining $D_E$
\textup{(}Definition~\ref{def:DE}\textup{)} and the generator
$\bar\gamma$ are fixed, the identity holds with a definite sign. Both
sides extend $\Q$-linearly in the symbol.
\end{theorem}

\begin{proof}
Brunault's intersection-product formula \cite[Prop.~19]{BrunaultSiegel}
states: for $x=\sum_i\{f_i,g_i\}\in K_2(E)\otimes\Q$ and
\emph{every} class $\gamma\in H_1(E(\C),\Z)$,
\begin{equation}\label{eq:brunault}
  r(x)[\gamma]\;=\;-\,(E^+\cdot\gamma)\;D_E\bigl(\lambda(x)\bigr),
  \quad
  \lambda(x)=\sum_i\sum_{p,q}\operatorname{ord}_p(f_i)\,
  \operatorname{ord}_q(g_i)\,\cl(p-q),
\end{equation}
where $E^+$ is the generator of $H_1(E,\Z)^+$ attached to the chosen
orientation of $E(\R)$ and $D_E$ is the elliptic dilogarithm of that
orientation. The divisor $\lambda(x)$ is term by term the diamond
product $\sum_i(f_i)\diam(g_i)$. By Lemma~\ref{lem:pairing}, a lift
$\gamma$ of the quotient generator $\bar\gamma$ satisfies
$E^+\cdot\gamma=\pm1$, and by Section~\ref{sec:setup} the left side of
\eqref{eq:brunault} depends only on the image of $\gamma$ in $Q(E)$.
Substituting gives \eqref{eq:main}. The $\Q$-linearity in the symbol is
the validity of \eqref{eq:brunault} for arbitrary decompositions
$x=\sum_i\{f_i,g_i\}$; in particular the right side is independent of
the decomposition, since the left side is.
\end{proof}

\begin{corollary}[The subgroup version, with the correction
factor]\label{cor:subgroup}
Let $\gamma^-$ be a generator of the anti-invariant subgroup
$H_1(E,\Z)^-$. Then, with $\iota(E)$ as in \eqref{eq:iota},
\begin{equation}\label{eq:subgroup}
  \int_{\gamma^-}\eta(f,g)
  \;=\;\pm\,\iota(E)\cdot D_E\bigl((f)\diam(g)\bigr).
\end{equation}
In particular, after matching orientations, the coefficient in
\eqref{eq:transmitted} is correct when $\Delta>0$; when $\Delta<0$
the magnitudes of its two sides differ by a factor of exactly $2$
for every symbol with
$D_E\bigl((f)\diam(g)\bigr)\neq0$.
\end{corollary}

\begin{proof}
By Lemma~\ref{lem:index}, $[\gamma^-]=\pm\iota(E)\,[\bar\gamma]$ in
$Q(E)$; apply Theorem~\ref{thm:main}.
\end{proof}

\begin{remark}\label{rem:17}
For $\Delta>0$, the subgroup and quotient generators have the same
normalization up to sign. This covers the conductor-$17$ computation
of \cite{LalinRam}, whose curve 17a4 has minimal discriminant $+17$.
For $\Delta<0$, Proposition~\ref{prop:11a3} gives an illustration by
proved values (with an external input from \cite{Paper1}), and
Section~\ref{sec:numerics} supplies a numerical illustration on 20.a4.
Touafek's printed formula on the latter curve is compatible with
the corrected normalization, because it includes the factor $\pi$
discussed in Remark~\ref{rem:touafek}.
\end{remark}

%=====================================================================
\section{An illustration and a comparison of conventions}\label{sec:proved}

\begin{proposition}[11.a3: comparison of proved values]\label{prop:11a3}
Let $E$ be the minimal model
$y_W^2+y_W=x_W^3-x_W^2$ of \textup{11.a3}, with $A=(0,0)$ and
origin $O$ at infinity. Let $x,y$ denote, instead, the coordinates
of $S_0:y^2+(x^2+1)y+x^3=0$, viewed as rational functions on $E$
through the birational identification in \cite{Paper1}.
Choose the real orientation used for the special values in
\cite[Th\'eor\`eme~8 and Corollaire~101, equation~(3.151)]{Brunault},
and write
\[
  b_{11}=L'(E,0)=\frac{11}{4\pi^2}L(E,2)>0.
\]
Here the equality follows from the functional equation (the curve has
root number $+1$), and positivity follows from the absolutely
convergent Euler product at $s=2$. Then
\begin{enumerate}
  \item For a generator $\gamma^-$ of $H_1(E,\Z)^-$,
  $\displaystyle\int_{\gamma^-}\eta(x,y)=\pm2\pi b_{11}$.
  This is the regulator-anchor theorem of \cite{Paper1}, used here
  as an external proved input; its Siegel-unit computation does
  not use the diamond normalization being compared here.
  \item $D_E\bigl((x)\diam(y)\bigr)=\frac52D_E(A)=\pi b_{11}$.
  Indeed, the divisor table from \cite{Paper1} is
  \[
    (x)=[A]+[2A]-[O]-[3A],\qquad
    (y)=3[2A]-2[3A]-[O].
  \]
  It gives $(x)\diam(y)=6(O)-5(A)+5(2A)$ in
  $\Z[E(\C)]^-$.
  The evaluations in
  \cite[Corollaire~101, equation~(3.151)]{Brunault} are
  \[
    D_E(2A)=\frac32D_E(A),\qquad
    L(E,2)=\frac{10\pi}{11}D_E(A).
  \]
  Combining the second with the functional equation gives
  $D_E(A)=\frac{2\pi}{5}b_{11}$ in the same orientation.
  The exotic relation is due to Bertin \cite{Bertin}.
\end{enumerate}
Thus the two evaluations realize the factor $\iota(E)=2$ of
Corollary~\ref{cor:subgroup}. This is a comparison using cited results,
not a new proof of the conductor-$11$ split-integral identity.
The proof of Theorem~\ref{thm:main} is independent of this example.
\end{proposition}

\begin{remark}[Touafek's printed normalization]\label{rem:touafek}
To avoid conflating conventions, write $r_T$ for the scalar regulator
in \cite[p.~119, Definition~2]{Touafek}:
\[
  r_T(\{f,g\})=\frac1{2\pi}\int_{\gamma^-}\eta(f,g),
  \qquad \gamma^-\ \text{a generator of }H_1(E,\Z)^-.
\]
His Theorem~1 and both identities of Proposition~1
\cite[pp.~121--122]{Touafek} have a factor $\pi$ on the left:
\begin{equation}\label{eq:touafekprinted}
  \pi r_T(\{f,g\})=D_E\bigl((f)\diam(g)\bigr).
\end{equation}
Here the elliptic dilogarithm is the usual $q$-series, with no extra
scalar in its definition \cite[p.~120, Definition~3]{Touafek}.
By Corollary~\ref{cor:subgroup}, the general comparison in these
conventions is instead
\begin{equation}\label{eq:touafekcorrect}
  \pi r_T(\{f,g\})
     =\pm\frac{\iota(E)}2D_E\bigl((f)\diam(g)\bigr).
\end{equation}
Consequently \eqref{eq:touafekprinted} is compatible with the
subgroup normalization on a curve with $\Delta<0$, after matching
orientations. This applies to the 20.a4 identity in his
Proposition~1:
\[
  \pi r_T(\{X_2,Y_2\})
     =6D_{E_2}(P_2)+6D_{E_2}(2P_2).
\]
For $\Delta>0$, the coefficient in
\eqref{eq:touafekcorrect} is $1/2$, up to sign, not $1$; the
difference in magnitudes is substantive whenever the diamond value
is nonzero. This is a comparison of normalizations, not an
independent nonvanishing claim about every symbol in that paper.

The earlier draft of this note omitted the printed $\pi$ from
Proposition~1 and consequently attributed an additional analytic
mismatch to it. That reading, and the proposed \emph{two-part
reinterpretation} based on it, are incorrect and are withdrawn.
Neither removing the factor $(2\pi)^{-1}$ from $r_T$ nor replacing
its subgroup generator is required for the $\Delta<0$ identity.
A full analysis of Touafek's isogeny argument in Theorem~2 and of
the derivation of his final exotic relation is not undertaken here.
In particular, we make no claim here that our normalization
comparison proves or disproves either of those statements.
\end{remark}

%=====================================================================
\section{A high-precision numerical instance on the curve
20.a4}\label{sec:numerics}

Our numerical example is the curve and symbol in Touafek's
Proposition~1 \cite{Touafek}. It illustrates the distinction between
the coefficient-$1$ subgroup reading \eqref{eq:transmitted} and
Touafek's differently normalized, compatible formula:
\[
  E_2:\;y^2+2xy+2y=(x-1)^3,
  \qquad\text{minimal model }y^2=x^3+x^2+4x+4.
\]
This is LMFDB \textup{20.a4} = Cremona \textup{20a1} \cite{LMFDB},
with $\Delta=-6400<0$. For the $6$-torsion point $P_2=(5,4)$ in
Touafek's model, his printed statement is
$\pi r_T(\{X_2,Y_2\})=6D_E(P_2)+6D_E(2P_2)$, not the same
identity with $r_T$ alone on the left.
We verified in PARI/GP \cite{PARI} that $P_2$ has order $6$,
that the divisor computations
$(X_2)=2[3P_2]-2[O]$, $(Y_2)=3[2P_2]-3[O]$ are correct, and
that $2\omega_2-\omega_1$ is the purely imaginary anti-invariant
lattice direction, as in Lemma~\ref{lem:basis}(2).

We then computed, at $80$ digits, the two regulator integrals of
$\eta(X_2,Y_2)$ --- over a lift of the \emph{quotient} generator and
over the \emph{subgroup} generator --- by an unwrapped-argument
midpoint scheme on loops shifted by a generic
$u_0=\omega_1/7+\omega_2/13$, with Richardson extrapolation ($p=2$)
and honest residuals (ancillary files \texttt{touafek\_20b.gp} and
\texttt{touafek\_20b\_log.txt}), and independently evaluated the
diamond side by the $q$-series \eqref{eq:DE} ($|q|=0.0797$, summed
symmetrically until the next paired term dropped below $10^{-75}$). The result:
\begin{align*}
  \text{quotient integral} &= -5.02110876405486(0)\ldots,\\
  6\bigl(D_E(P_2)+D_E(2P_2)\bigr) &= -5.0211087640548395\ldots,\\
  \text{residual} &= -2.1\times10^{-14},\\
  \text{subgroup integral}/\text{quotient integral}
    &= 1.9999999999999771\ldots
\end{align*}
The quotient integral agrees with the diamond value (factor $1$) to
$2.1\times10^{-14}$ --- consistent with the expected $N^{-4}$
remainder of the extrapolation,
well within the method's demonstrated error $e(N)\approx2.98\,N^{-2}$
before extrapolation --- and the subgroup integral agrees with twice
the quotient integral to $1.2\times10^{-13}$; the transmitted
reading \eqref{eq:transmitted} fails by a full diamond unit.
In contrast, $\pi r_T$ is half the subgroup integral and therefore
agrees numerically with the diamond value, with the chosen
orientations. The computation illustrates
Theorem~\ref{thm:main} and Corollary~\ref{cor:subgroup}; it does not
provide a counterexample to Touafek's printed Proposition~1 on
20.a4.

\begin{remark}[Status of the numerics]\label{rem:status}
The two sides of this section have different error profiles. The
diamond side is the $q$-series \eqref{eq:DE} with $|q|=0.0797$: the
symmetric summation was stopped when the next paired term had size
below $10^{-75}$, and the displayed value was stable under increased
precision and truncation; it is not an interval-certified value. The
integral side is a midpoint rule on the
unwrapped argument: its error before extrapolation is
$e(N)\approx 2.98\,N^{-2}$ (measured by refinement), and one
Richardson step ($p=2$) eliminates the $N^{-2}$ term; the residual
$-2.1\times10^{-14}$ against the $q$-series value is consistent with
an $N^{-4}$ remainder. The digits printed in this section are
those stable under $N\mapsto 2N$ refinement followed by extrapolation;
they should not be read as $20$ certified digits. Altogether the
computations are high-precision numerical
confirmations (agreement to $\ge13$ digits between two independent
evaluations), not interval-arithmetic certificates; they play no role
in the proofs of Section~\ref{sec:corrected}, which rest on the
published, proved \cite[Prop.~19]{BrunaultSiegel}. The scripts are
self-contained and reproducible.
\end{remark}

%=====================================================================
\section{The transmission of the statement}\label{sec:transmission}

We record, neutrally, what we have been able to establish about the
provenance of the wording in \eqref{eq:transmitted}.

(1) \emph{Bloch's 1978 ICM report} \cite{BlochICM} contains the diamond
construction in its original form: the sum
$\sum_{i,j}m_in_jD_q(a_i/b_j)$ is shown to be well-defined and to induce
$K_2(\C(E))\to\R$ directly (the full $R_q$ maps to $\C$; its real part
is the diamond sum). To the best of our reading, \emph{no loop $\gamma$
appears}: in this form there is no lattice to choose and nothing to
misnormalize. The same is true of the Eisenstein--Kronecker treatment
of $L(E,2)$ by Goncharov--Levin \cite{GoncharovLevin}, which proceeds
through the diamond picture without invoking a generator of
$H_1(E,\Z)^-$.

(2) \emph{The cohomological form} $r:K_2(E)\to H^1(E,\R)$,
$r(\{f,g\})[\gamma]=\int_\gamma\eta(f,g)$, makes the regulator a linear
form on homology. Since it vanishes on $H_1(E,\Z)^+$, its natural
integral lattice is $Q(E)$ --- this is the reading proved in
Theorem~\ref{thm:main} and, in intersection-product form, already in
\cite[Prop.~19]{BrunaultSiegel}. Brunault himself applies it on curves
with $\Delta<0$, printing $E^+\cdot\gamma=-2$ for a generator of
$H_1(E,\Z)^-$ \cite[\S6.1]{BrunaultSiegel}: the factor $2$ was in
print. What we have not found in print is the observation that the
subgroup wording of \eqref{eq:transmitted}, as in
\cite[Thm.~6]{LalinRam} and the papers following it, is incompatible
with that normalization when $\Delta<0$.

(3) \emph{Analytic prefactors must be read together with the cycle
normalization.} Deninger \cite{Deninger} does not use the diamond
notation: he writes
$m(P)=-\langle r_D[t_1,t_2]_E,\gamma\rangle$, with a cycle dual to a
generator of $H^1(E/\R,\Q(1))$ \cite[p.~274]{Deninger}.
The formula $m(P)=(2\pi)^{-1}r(\{x,y\})[\gamma]$ is the later
reformulation in \cite[eq.~(11)]{LalinRam}.
For Touafek, the simultaneous presence of $(2\pi)^{-1}$ in
Definition~2 and of $\pi r_T$ in Theorem~1 and Proposition~1 gives
$\pi r_T=\frac12\int_{\gamma^-}\eta$. On a curve with $\Delta<0$
this has the same magnitude as the diamond value by
Corollary~\ref{cor:subgroup}. There is no additional factor of $\pi$
to remove. Thus the incorrect coefficient-$1$ subgroup formula
\eqref{eq:transmitted} must not be confused with Touafek's printed
formula \eqref{eq:touafekprinted}. Our comparison does not settle
his separate isogeny argument; see Remark~\ref{rem:touafek}.

(4) \emph{We have not adjudicated Bloch's book.} The CRM monograph
\cite{Bloch} itself was not available to us; whether the subgroup
wording with factor $1$ originates there or in a subsequent
formulation is not decided here. We do not claim to have established
a complete chain of textual transmission. Nothing in
this note depends on the answer.

%=====================================================================
\section*{Statements and Declarations}

\paragraph{Funding} The author did not receive support from any
organization for the submitted work.

\paragraph{Competing Interests} The author has no relevant financial or
non-financial interests to disclose.

\paragraph{Data Availability}
The verification script and its numerical run log
(\path{touafek_20b.gp}, \path{touafek_20b_log.txt}) are provided in
the \path{anc/} directory of the source distribution.
The script comments have been corrected to include the printed
$\pi$ in Touafek's Proposition~1; the numerical calculation concerns
bare integrals and is unchanged.
The conductor-$11$ example uses the external results and companion
artifacts of \cite{Paper1}; they are not reproved in this note.

\paragraph{Use of AI tools} Documented in the Declaration on the use of
AI tools preceding the Introduction.

%=====================================================================

\end{document}